\documentclass[12pt]{amsart}
\usepackage{amsthm}
\usepackage{amssymb,amsmath,mathtools}
\newtheorem{proposition}{Proposition}[section]

\theoremstyle{definition}

\newtheorem{remark}[proposition]{Remark}

\numberwithin{equation}{section}
\def \e {\mathsf{e}}
\def \q {\boldsymbol{q}}
\def \x {\boldsymbol{x}}

\def\R {\mathbb{R}}

\def \au {\rm}
\def \ti {\it}
\def \jou {\rm}
\def \bk {\it}
\def \no#1#2#3 {{\bf #1} (#3), #2.}
\def \eds#1#2#3 {#1, #2, #3.}
\title[Heat conduction laws]
{A hierarchy of heat conduction laws}

\author[F. Dell'Oro and V. Pata]
{Filippo Dell'Oro and Vittorino Pata}

\address{Politecnico di Milano - Dipartimento di Matematica
\newline\indent
Via Bonardi 9, 20133 Milano, Italy}
\email{filippo.delloro@polimi.it {\rm (F. Dell'Oro)}}
\email{vittorino.pata@polimi.it {\rm (V. Pata)}}


\begin{document}

\begin{abstract}
The purpose of this work is to produce a family of equations
describing the evolution of the temperature in a rigid heat conductor.
This is obtained by means of successive approximations
of the Fourier law, via memory relaxations and integral perturbations.
\end{abstract}

\maketitle

\section{Introduction}

\noindent
Consider a homogenous isotropic rigid heat conductor occupying
a domain $\Omega\subset\R^N$
and denote by
$$u=u(\x,t):\Omega\times\R^+\to\R$$
its \emph{temperature variation} at each point and time
from an equilibrium reference value $\theta_0$, that is,
$$
u(\x,t) = {{\theta(\x,t) - \theta_0}\over {\theta_0}},
$$
where $\theta$ is the absolute temperature.
In absence of heat energy generation, from external or internal sources,
the thermal evolution of the body
is governed by the energy balance equation
$$
\partial_t \e + {\rm div}\, \q=0.
$$
Here,
$$\e=\e(\x,t):\Omega\times\R^+\to\R$$
is
the \emph{internal energy} of the conductor, while
$$\q=\q(\x,t):\Omega\times\R^+\to \R^N$$
is the \emph{heat flux vector}.
Within the assumption of small variations of the temperature and its gradient,
the rate of change of the internal energy in the material
is proportional to the rate of change of its temperature (see \cite{GUP}).
Namely,
$$\partial_t \e = c \partial_t u,$$
where the proportionality constant
$c>0$ is the \emph{specific heat}.
Accordingly, setting from now on $c=1$,
the energy balance equation takes the form
\begin{equation}
\label{BASE}
\partial_t u + {\rm div}\, \q=0.
\end{equation}
In order to derive an equation ruling the evolution of the
relative temperature, one more ingredient is needed:
the so-called constitutive law for the heat flux,
allowing to write $\q$ in terms of $u$. In fact,
the choice of
the constitutive law is
what really determines a model of heat conduction.
At the same time, being a purely heuristic
interpretation of a certain  physical phenomenon, it
may reflect different individual perceptions
of reality, or even philosophical beliefs.

\smallskip
The first constitutive law for the heat flux appeared in the literature
is $\q=-\kappa\nabla u$, with $\kappa>0$, established in 1882 by Fourier 
in his celebrated work \emph{Th\'eorie analytique de la chaleur}~\cite{FOU}.
Since then, as shown in the next section, a number of variants
of this law have been proposed.
In particular, in 1995, Tzou suggested in the famous paper~\cite{TZO} a modification
of the Fourier law where both $\q$ and the gradient of temperature 
present a delay, namely,
\begin{equation}
\label{UNO}
\q(\x,t+\tau_{\q})=-\kappa\nabla u (\x,t+\tau_u),
\end{equation}
where $\tau_{\q},\tau_u>0$ are delay parameters.
Some years later, Tzou’s theory has been modified via 
the introduction of a constitutive equation of form \cite{RCH}
\begin{equation}
\label{DUE}
\q(\x,t+\tau_{\q})=-\kappa\nabla u (\x,t+\tau_u)-\kappa_1\nabla u_1 (\x,t+\tau_{u_1}),
\end{equation}
where $\tau_{u_1}$ is another delay parameter which is also assumed to be positive.
Here, $u_1$ is the so-called thermal displacement satisfying $\partial_t u_1=u$,
and $\kappa_1>0$ is a constant which is typical in the type II and III thermoelastic theories.
Heat conduction laws of the above kind are known as
\emph{phase lag models}; more precisely, dual phase lag if~\eqref{UNO}, three phase lag if~\eqref{DUE}.
Since the delay parameters are assumed to be very small,
Quintanilla and Racke~\cite{QUR} (but see also the textbook~\cite{STR})
proposed to take
the formal Taylor expansions of~\eqref{UNO} and~\eqref{DUE}
stopping at different orders, that is
(omitting the dependence on the space variable $\x$),
\begin{align*}
&\q(t)+\tau_{\q}\partial_t\q(t)+\frac12\tau_{\q}^2
\partial_{tt}\q(t)+\ldots\\
&=-\kappa\nabla u (t)
-\tau_u\kappa\nabla \partial_t u (t)
-\frac12\tau_u^2\kappa \nabla \partial_{tt} u (t)+\ldots,
\end{align*}
and
\begin{align*}
&\q(t)+\tau_{\q}\partial_t\q(t)+\frac12\tau_{\q}^2
\partial_{tt}\q(t)+\ldots\\
&=-\kappa\nabla u (t)
-\tau_u\kappa\nabla \partial_{t} u (t)
-\frac12\tau_u^2\kappa\nabla\partial_{tt} u (t)+\ldots\\
&\quad -\kappa_1\nabla u_1 (t)
-\tau_{u_1}\kappa_1\nabla \partial_{t} u_1 (t)
-\frac12\tau_{u_1}^2\kappa_1\nabla \partial_{tt} u_1 (t)+\ldots.
\end{align*}
Therefore, the sum 
$$\eqref{BASE}+\tau_{\q}\partial_t \eqref{BASE}+\frac12\tau_{\q}^2
\partial_{tt}\eqref{BASE}+\ldots
$$
gives rise to a hierarchy
of evolution equations in the variable $u$ and $u_1$, respectively.

\smallskip
In the present work, taking inspiration from this idea,
we propose an alternative path, which allows in particular to include
into the picture some heat laws containing memory terms
described by convolution integrals,
hence having a non local character in time.
The starting point is again the Fourier law. At each step, we find a heat law involving 
$\nabla u_n$ along with its time derivatives
up to a certain order,
where $\partial_t^n u_n=u$. Then, we perform a memory relaxation
of $\nabla u_n$ via a convolution integral, and we introduce a further integral
perturbation by adding the term $-\kappa_{n+1}\nabla u_{n+1}$, with $\kappa_{n+1}>0$ small.
When the convolution kernel is the negative exponential,
we can show that what we obtain is nothing but the analogue of the
law of the previous step,
now expressed in terms of $u_{n+1}$, 
and where the order of the time derivatives of the heat flux has increased
by one.
We depict our procedure in the following flow chart:

\smallskip
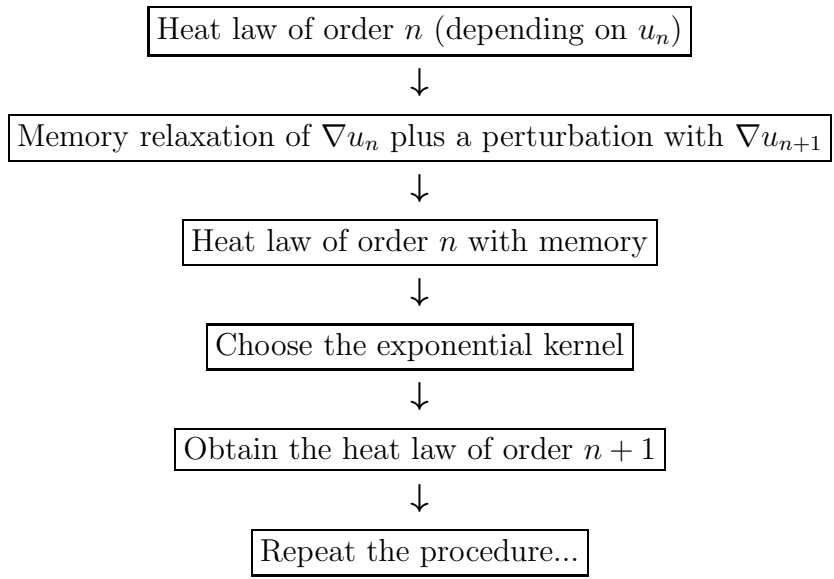
\begin{gather*}
\fbox{\text{Heat law of order $n$ (depending on $u_n$)}}\\
\boldsymbol{\downarrow}\\
\fbox{\text{Memory relaxation of $\nabla u_n$ plus a perturbation
with $\nabla u_{n+1}$}}\\
\boldsymbol{\downarrow}\\
\fbox{\text{Heat law of order $n$ with memory}}\\
\boldsymbol{\downarrow}\\
\fbox{\text{Choose the exponential kernel}}\\
\boldsymbol{\downarrow}\\
\fbox{\text{Obtain the heat law of order $n+1$}}\\
\boldsymbol{\downarrow}\\
\fbox{\text{Repeat the procedure...}}
\end{gather*}

\medskip
\noindent
With this method, we produce a hierarchy of evolution equations,
along with their memory counterparts. In fact, most of our equations without memory
are the same of those found in~\cite{QUR}. However, it is important to point out
that, contrary to~\cite{QUR} where the equations are 
at any step in the variables $u$ or $u_1$, in our case
the equation at the $n^{\rm th}$-step is for the variable $u_n$.

\section{The Fourier Law and its Relaxations}

\subsection{The Fourier law}
For the classical Fourier constitutive law~\cite{FOU}
\begin{equation}
\label{Fourier}
\q=-\kappa\nabla u,
\end{equation}
where $\kappa>0$ is the \emph{instantaneous conductivity},
we deduce from~\eqref{BASE} the well-known
heat equation
\begin{equation}
\label{FourierE}
\partial_t u-\kappa\Delta u=0.
\end{equation}
This equation, being parabolic, exhibits a physically unpleasant feature:\
a disturbance at any point is felt instantaneously
everywhere. In other words, the speed
propagation of thermal
signals is infinite.
This phenomenon,
sometimes referred to as the \emph{paradox
of heat conduction} (see, e.g., \cite{CJ,FIC}),
is not expected, nor observed, in real conductors.
Therefore, several attempts have been made through the years
in order to introduce some hyperbolicity in the mathematical modeling
of heat conduction.
Starting from~\eqref{Fourier}, the idea is then to produce more
refined heat laws by means of perturbative arguments,
or by applying some kind of relaxation
to the variables in play.

\subsection{The Maxwell-Cattaneo law}
A first possibility is adopting the
renowned Maxwell-Cattaneo law~\cite{CAT},
namely, the differential perturbation
of~\eqref{Fourier}
\begin{equation}
\label{MC}
\q+\varepsilon \partial_t\q=-\kappa\nabla u,
\end{equation}
where $0<\varepsilon\ll \kappa$.
In which case, the sum
$\eqref{BASE}+\varepsilon\partial_t\eqref{BASE}$
entails
the weakly damped wave equation
\begin{equation}
\label{WDWE}
\varepsilon \partial_{tt}u+\partial_t u-\kappa\Delta u=0,
\end{equation}
otherwise called the telegrapher's equation, widely employed in the description
of several physical phenomena.

\subsection{The Gurtin-Pipkin and the Coleman-Gurtin laws}
Another strategy is
to relax the term (or part of the term) $\nabla u$ appearing
in the Fourier law~\eqref{Fourier} by means of a time-convolution
against a
bounded convex summable function $g:\R^+\to\R^+$, 
usually called \emph{memory kernel}, of total mass
$$\int_0^\infty g(s)ds=1.$$
Precisely, for $\omega\in[0,1)$, we take (omitting the dependence on $t$ outside the integral)
\begin{equation}
\label{GPCG}
\q=-\omega\kappa\nabla u- (1-\omega)\kappa \int_0^\infty g(s)\nabla u(t-s)ds =0,
\end{equation}
where the initial past history of the temperature $u(\x,t)_{|t\leq 0}$
is supposed to be known, and it is regarded as an initial datum of the problem.
The above is usually called the Gurtin-Pipkin law~\cite{GUP} when $\omega=0$,
and the Coleman-Gurtin law~\cite{COG} when $\omega\in(0,1)$.
Plugging~\eqref{GPCG} into \eqref{BASE}
we end up with the integro-differential equation
\begin{equation}
\label{GPCGE}
\partial_t u-\omega\kappa\Delta u-(1-\omega)\kappa \, \int_0^\infty g(s)\Delta u(t-s)ds=0,
\end{equation}
which exhibits a nonlocal (in time) character, due to the memory term.
The role of $\omega$ is most peculiar: when $\omega>0$, the presence
of the instantaneous Laplacian of the temperature has a regularizing effect, introducing parabolicity in the model.
On the contrary, the equation for $\omega=0$ is purely hyperbolic (see \cite{FGR,GP,GUP}).
In fact, this regularizing feature will occur in all the heat models considered hereafter, whenever $\omega>0$.
Finally, let us observe that
in the formal limiting case where $\kappa$ equals the Dirac mass at zero,
or when $\omega=1$,
we recover \eqref{FourierE}. From the physical viewpoint, this means that
\eqref{GPCGE} is close to \eqref{FourierE} when the memory kernel is concentrated, that is,
when the system keeps a very short memory of past effects.

\begin{remark}
Choosing $u(\x,t)_{|t\leq 0}=0$, that is, null initial past history of the temperature,
one obtains the Volterra version of~\eqref{GPCGE}
$$\partial_t u-\omega\kappa\Delta u-(1-\omega)\kappa \, \int_0^t g(s)\Delta u(t-s)ds=0.$$
\end{remark}

\subsection{The exponential kernel}
\label{subexp}
We now consider the heat law~\eqref{GPCG} for the
exponential kernel
\begin{equation}
\label{gexp}
g(s)=\frac1\varepsilon\, e^{-\frac{s}\varepsilon},\quad\varepsilon>0.
\end{equation}
Then, observing that
\begin{align}
\label{dercon}
\partial_t \int_0^\infty g(s)\nabla u(t-s)ds
&=\frac{1}{\varepsilon}\partial_t \,\bigg(e^{-\frac{t}\varepsilon}\int_{-\infty}^t e^{\frac{s}\varepsilon}\nabla u(s)ds\bigg)\\
&=\frac{1}{\varepsilon} \nabla u(t)-\frac{1}{\varepsilon}\int_0^\infty g(s)\nabla u(t-s)ds,\notag
\end{align}
taking the sum $\eqref{GPCG}+\varepsilon\partial_t\eqref{GPCG}$
we end up with
\begin{equation}
\label{MC2}
\q+\varepsilon \partial_t\q=-\kappa\nabla u-\varepsilon\omega\kappa\nabla \partial_t u,
\end{equation}
which contains the Maxwell-Cattaneo law~\eqref{MC} as a particular instance,
corresponding to $\omega=0$.
In particular, this tells that
a differential perturbation of the Fourier law
is just a memory relaxation with the exponential kernel.
Again, by $\eqref{BASE}+\varepsilon\partial_t\eqref{BASE}$, we
obtain the evolution equation for the temperature, that now reads
\begin{equation}
\label{WSDWE}
\varepsilon \partial_{tt}u+\partial_t u
-\varepsilon\omega\kappa\Delta\partial_t u-\kappa\Delta u =0.
\end{equation}
Clearly, when $\omega=0$ this is the weakly damped wave equation~\eqref{WDWE}.
On the other hand, when $\omega>0$ we have the regularizing extra term
$-\Delta\partial_t u$, which renders the equation no longer purely hyperbolic.
As a matter of fact, we can also
allow the limiting situation $\varepsilon=0$, corresponding to
the Dirac mass at zero. In which case, we boil down to the Fourier law,
and \eqref{WSDWE} becomes the heat equation~\eqref{FourierE}.

\begin{remark}
\label{remmoremmo}
Apparently, for the limit value $\omega=1$, this procedure might look flawed. Indeed,
from the one side if $\omega=1$ the starting law~\eqref{GPCG} is just
the Fourier law~\eqref{Fourier}. From the other side the procedure
above leads to~\eqref{MC2} with $\omega=1$, namely,
$$
\q+\varepsilon \partial_t\q=-\kappa\nabla u-\varepsilon\kappa\nabla \partial_t u.
$$
But, calling
$$\boldsymbol{\phi}=\q+\kappa\nabla u,$$
the latter equation can be written as
$$
\partial_t \boldsymbol{\phi}+\frac1\varepsilon \boldsymbol{\phi}=0,
$$
implying in turn that
$$\boldsymbol{\phi}(t)=\boldsymbol{\phi}(0)e^{-\frac{t}\varepsilon}.$$
Note that this is possible just because $\omega=1$. At this point,
all depends on the initial value $\boldsymbol{\phi}(0)$ which must be necessarily read
from the original equation~\eqref{GPCG}, clearly with $\omega=1$.
This yields $\boldsymbol{\phi}(0)=0$, hence $\boldsymbol{\phi}=0$ for all times,
recovering the Fourier law.
This is the reason why we always take $\omega\in[0,1)$.
\end{remark}

\subsection{Heat conduction of type III}
The theory of heat conduction of type III goes back to the works of
Green and Naghdi~\cite{GN1,GN2,GN3,GN4,GN5}, and it is based on the introduction of a further variable
$$u_1=u_1(\x,t):\Omega\times\R^+\to\R,$$
called \emph{thermal displacement}, defined (up to a constant in time function)
by the relation
$$\partial_t u_1=u.$$
The Fourier law~\eqref{Fourier}
is modified as
\begin{equation}
\label{GN}
\q=-\kappa\nabla u-\kappa_1\nabla u_1,
\end{equation}
where the constant $\kappa_1>0$ is the
\emph{conductivity rate}.
Since
$$u_1(\x,t)=u_1(\x,0)+\int_0^t u(\x,s)ds,$$
the above can be viewed as an integral perturbation
of~\eqref{Fourier}, which collapses into~\eqref{Fourier} in the limit case $\kappa_1=0$. From
the energy balance~\eqref{BASE}, we deduce the evolution equation for the variable $u_1$
\begin{equation}
\label{SDWE}
\partial_{tt}u_1 -\kappa\Delta\partial_t u_1-\kappa_1\Delta u_1=0.
\end{equation}
Equation~\eqref{SDWE}, known as the strongly damped wave equation,
has applications in several fields. In particular, it serves as a model of viscoelasticity,
more precisely, in the language of Dautray and Lions~\cite{DOLI},
of viscoelasticity with \emph{short memory}.
In that context, it goes by the name of
Kelvin-Voigt equation.

\section{Heat Law of Order 1}

\noindent
The idea now is to consider the modifications \eqref{GPCG} and \eqref{GN}
of the Fourier law at the same time. This leads to the
introduction of the following heat law, that within the language of the next section we might call
\emph{of order $0$ with memory}:
\begin{equation}
\label{0mem}
\q=-\omega\kappa\nabla u- (1-\omega)\kappa \int_0^\infty g(s)\nabla u(t-s)ds
-\kappa_1\nabla u_1=0.
\end{equation} From the energy balance~\eqref{BASE},
and recalling that $u=\partial_t u_1$,
we find the integro-differential equation in the variable $u_1$
\begin{equation}
\label{0memE}
\partial_{tt} u_1-\omega\kappa\Delta \partial_t u_1-(1-\omega)\kappa \, \int_0^\infty g(s)\Delta \partial_t u_1(t-s)ds-\kappa_1\Delta u_1=0.
\end{equation}
Equation \eqref{0memE}
for $\omega=0$ is the renowned equation of viscoelasticity
with \emph{long memory}~\cite{DOLI}, which is
purely hyperbolic, contrary to the its short memory counterpart, formally obtained
when $g$ is the Dirac mass at zero.
This equation has been first analyzed in a systematic way by Dafermos~\cite{DAF},
and studied by quite many authors thereafter. Instead,
when $\omega>0$ we have the strongly
damped wave equation with memory, where one can appreciate the regularization provided by the
term $-\Delta \partial_t u_1$, that confers the equation a partially parabolic character
(see~\cite{DPZ}).

\begin{remark}
Actually, introducing the differentiated memory kernel
$$\mu(s)=-\kappa g'(s),$$
by means of an integration by parts
the equation of viscoelasticity can be given the more familiar form
$$
\partial_{tt} u_1-\kappa_1\Delta u_1-\int_0^\infty \mu(s)
\big[\Delta u_1(t)-\Delta u_1(t-s)\big]ds=0.
$$
Indeed, the boundary terms of the integration by parts vanish within
reasonable regularity hypotheses on the solution (see, e.g., \cite{GPT}).
\end{remark}

The next step is repeating the scheme devised
in Subsection~\eqref{subexp}, by choosing the memory kernel $g$
of the form~\eqref{gexp}.
Accordingly, by taking the sum $\eqref{0mem}+\varepsilon\partial_t\eqref{0mem}$
and exploiting~\eqref{dercon},
we arrive at the heat law \emph{of order} $1$
\begin{equation}
\label{order1}
\q+\varepsilon\q_t=-\kappa_1\nabla u_1
-(\kappa+\varepsilon \kappa_1)\nabla \partial_t u_1-\varepsilon\omega\kappa\nabla \partial^2_t u_1,
\end{equation}
written in terms of the variable $u_1$ solely.
In turn, the energy balance~\eqref{BASE} produces the third-order
evolution equation
\begin{equation}
\label{MGT}
\varepsilon\partial_{ttt}u_1+\partial_{tt}u_1
-\varepsilon\omega\kappa\Delta \partial^2_t u_1
-(\kappa+\varepsilon \kappa_1)\Delta \partial_t u_1-\kappa_1\Delta u_1=0.
\end{equation}
When $\omega=0$, this is the Moore-Gibson-Thompson (MGT) equation,
arising in the modeling of wave propagation
in viscous thermally relaxing fluids \cite{MGI,THO}, whose
first derivation is due to Stokes~\cite{STO}.
The case $\omega>0$ yields a partially parabolic
regularization of the (fully hyperbolic) MGT equation.

\begin{remark}
As in Subsection~\eqref{subexp}, we could also allow $\varepsilon$ to be zero,
corresponding to a kernel $g$ collapsing into the Dirac mass at zero.
Then~\eqref{order1} reduces to the law~\eqref{GN} of heat conduction of type III.
Clearly, if also $\kappa_1=0$, we recover the Fourier law expressed in terms
of $\partial_t u_1$, which is nothing but $u$.
\end{remark}

\subsection*{On the MGT equation}
Before going any further, let us dwell for a while on the MGT equation,
which in the literature is usually written in the form
(here in the variable $v$)
\begin{equation}
\label{MGTr}
\partial_{ttt}v+a\partial_{tt}v-b\Delta \partial_t v-c\Delta v=0,
\end{equation}
where $a,b,c$ are generic positive constants.
Provided that $\Omega$ has a sufficiently regular boundary,
assuming the homogeneous Dirichlet boundary condition for $v$,
equation~\eqref{MGTr} is well known
to generate a linear
solution semigroup $S(t)$
on the natural weak energy space
$${\mathcal H}=H_0^1(\Omega)\times H_0^1(\Omega)\times L^2(\Omega).$$
Here, with standard notation, $L^2$ is the Lebesgue space of square-summable functions,
whereas $H_0^1$ is the Sobolev space of functions which are
square-summable along with their
first derivatives and vanish on the boundary of $\Omega$.
Notably, the asymptotic properties of $S(t)$ depend on a combination
of the physical parameters. Indeed, defining the
so-called stability number
$$
\varkappa=b-\frac{c}{a},
$$
it turns out that the semigroup is exponentially stable if and only if
$\varkappa>0$ (subcritical regime),  whereas in the supercritical regime
$\varkappa<0$ there are trajectories whose energy blows up exponentially fast~\cite{KLM,MMT}.
This fact somehow indicates
that the MGT equation is physically meaningful only when $\varkappa>0$.
If we want to write our equation~\eqref{MGT} in the form~\eqref{MGTr} for
arbitrary coefficient $a,b,c>0$, we have to choose the constants as follows:
$$\omega=0,\qquad\varepsilon=\frac1a,\qquad \kappa_1=\frac{c}a,\qquad
\kappa=\frac{\varkappa}a.
$$
Since $\kappa$ must be positive, it is apparent that such a procedure allows us to recover
the MGT equation~\eqref{MGTr} \emph{only} in the subcritical case. And this actually makes
perfectly sense, as what we produced is an equation describing the evolution of the temperature, whose semigroup is supposed
to exhibit a fast decay in absence of heat sources.
This pattern has been already noted in~\cite{DPMGT}, where the MGT equation is viewed
as a particular instance of the equation of viscoelasticity.

\section{Heat Laws of Order $n$}

\noindent
Renaming the relative temperature $u$ by $u_0$, for $n\geq 1$ we introduce the sequence of functions $u_n$,
each one being the antiderivative of $u_{n-1}$, hence satisfying the relations
$$\partial_t u_n=u_{n-1}.$$
Our purpose, as mentioned in the Introduction, is to generate a hierarchy of heat laws
involving the function $u_n$, for every fixed $n$.
To this end,
let us first rename the constants $\varepsilon$, $\omega$ and $\kappa$ appearing
in~\eqref{order1} by $\varepsilon_1$, $\omega_1$ and $\kappa_0$, respectively.
Accordingly, the law of order~1 becomes
\begin{equation}
\label{order1new}
\q+\varepsilon_1\q_t=-\kappa_1\nabla u_1
-(\varepsilon_1 \kappa_1+\kappa_0)\nabla \partial_t u_1-\varepsilon_1\omega_1\kappa_0\nabla \partial^2_t u_1.
\end{equation}
We also assume that, for every $n\geq 1$, we are given the sequences
$$\varepsilon_n\geq 0 \qquad\text{and}\qquad\omega_n\in[0,1),$$
and, for every $n\geq 0$,
$$\kappa_n> 0.$$
Actually, we can allow $\kappa_{m}=0$ for some $m\geq 1$, but in that case
the sequences $\varepsilon_n,\omega_n,\kappa_n$ are finite and they stop at $n=m$.

Then, for every $n\geq 1$
(and whenever $\kappa_{n-1}>0$),
we define the heat law \emph{of order} $n$ by
\begin{equation}
\label{ordern}
\q+\sum_{i=1}^n  \alpha_n^i\partial_t^i\q=-\kappa_n\nabla u_n
-\sum_{i=1}^{2n}\beta_n^i\nabla \partial_t^i u_n.
\end{equation}
For $n=1$, in compliance with~\eqref{order1new},
the (nonnegative) coefficients appearing in the formula above are
$$\alpha_1^1=\varepsilon_1,\qquad
\beta_1^1=\varepsilon_1 \kappa_1+\kappa_0,\qquad
\beta_1^2=\varepsilon_1\omega_1\kappa_0.
$$
Instead, for $n\geq 2$ they are determined by the recurrence relations
\begin{equation}
\label{alfa}
\alpha_n^i=
\begin{cases}
\varepsilon_n+\alpha_{n-1}^1 &\text{ if }i=1,\\
\varepsilon_n\alpha_{n-1}^{i-1}+\alpha_{n-1}^i &\text{ if }n\geq 3 \text{ and }i=2,\ldots,n-1,\\
\varepsilon_n\alpha_{n-1}^{n-1} &\text{ if }i=n,
\end{cases}
\end{equation}
and
\begin{equation}
\label{beta}
\beta_n^i=
\begin{cases}
\varepsilon_n\kappa_n+\kappa_{n-1} &\text{ if }i=1,\\
\varepsilon_n\omega_n\kappa_{n-1}+\beta_{n-1}^1 &\text{ if }i=2,\\
\varepsilon_n\beta_{n-1}^{i-2}+\beta_{n-1}^{i-1} &\text{ if }i=3,\ldots,2n-1,\\
\varepsilon_n\beta_{n-1}^{2n-2} &\text{ if }i=2n.
\end{cases}
\end{equation}

\begin{remark}
\label{remalfa}
The coefficients $\alpha_n^i$ have the explicit form
$$
\alpha_n^i=\frac{1}{i!}\sum_{k_1,\ldots,k_i=1}^n
\varepsilon_{k_1}\cdots \varepsilon_{k_i}\delta_{k_1,\ldots, k_i},
$$
where
$$\delta_{k_1,\ldots, k_i}=
\begin{cases}
1 & \text{if $k_1,\ldots,k_i$ are distinct integers},\\
0 & \text{otherwise},
\end{cases}
$$
is the generalized delta Kronecker.
The relevant information provided by this formula is that if $\alpha_n^i=0$, then at least
$n-i+1$ among the constants $\varepsilon_1,\ldots,\varepsilon_n$
vanish.
Accordingly,
if $\alpha_n^{i}=0$ for some $i$, then $\alpha_n^k=0$ for all $k\geq i$.
\end{remark}

\begin{remark}
Concerning the coefficients $\beta_n^i$, one can check the following facts:
\begin{itemize}
\smallskip
\item[(i)] $\beta_n^i> 0$ for all $i=1,\ldots,n$.
\smallskip
\item[(ii)] $\beta_n^{2n}>0$ if and only if $\omega_1>0$ and $\varepsilon_1,\ldots,\varepsilon_n>0$.
\smallskip
\item[(iii)] Let $k=0,\ldots,n-1$ be fixed, and assume that $\alpha_{n}^{n-k}=0$. In particular, 
Remark~\ref{remalfa} implies that $\alpha_n^n=\alpha_n^{n-1}=\ldots=\alpha_{n}^{n-k}=0$.
Then, via (actually nontrivial) combinatorial arguments, it follows that
$$\beta_n^{2n}=\ldots\beta_n^{2n-k}=0.$$
\end{itemize}
\end{remark}

At this point, we distinguish two cases.

\medskip
\noindent
$\bullet$ If $\kappa_n=0$, we rewrite \eqref{ordern} in terms of the variable $u_{n-1}$,
so obtaining
$$
\q+\sum_{i=1}^n  \alpha_n^i\partial_t^i\q=
-\sum_{i=0}^{2n-1}\beta_n^{i+1}\nabla \partial_t^{i} u_{n-1}.
$$
Hence, taking the sum $\eqref{BASE}+\sum_{i=1}^n  \alpha_n^i\partial_t^i\eqref{BASE}$,
we find the evolution equation in the variable $u_{n-1}$
\begin{equation}
\label{nmeno}
\sum_{i=1}^n  \alpha_n^i\partial_t^{i+n}u_{n-1}+\partial_t^{n}u_{n-1}
-\sum_{i=0}^{2n-1}\beta_n^{i+1}\Delta \partial_t^{i} u_{n-1}=0.
\end{equation}
In particular, the last term $\beta_n^{2n}\nabla \partial_t^{2n-1} u_{n-1}$,
if present, is responsible for a regularizing effect. On the contrary, the character of
\eqref{nmeno} is purely hyperbolic when all the coefficients $\alpha_n^i$ differ from zero and
$\beta_n^{2n}=0$.
As a particular case of \eqref{nmeno}, we like to mention the case $n=2$ with $\beta_2^4=0$ (i.e., in absence of regularization), and all the other coefficients strictly positive, namely,
$$
\alpha_2^2\partial_{tttt}u_{1} +\alpha_2^1\partial_{ttt}u_{1}+\partial_{tt}u_{1}
-\beta_2^{3}\Delta \partial_{tt} u_{1}-\beta_2^{2}\Delta \partial_{t} u_{1}
-\beta_2^{1}\Delta u_{1}=0.
$$
This is a fourth-order equation of
MGT type which has been extensively studied in~\cite{MGT4}.

\medskip
\noindent
$\bullet$ If $\kappa_n\neq 0$,
we find the evolution equation in the variable $u_{n}$
\begin{equation}
\label{npiu}
\sum_{i=1}^n  \alpha_n^i\partial_t^{i+n+1}u_{n}+\partial_t^{n+1}u_{n}
-\sum_{i=1}^{2n}\beta_n^i\Delta \partial_t^i u_n-\kappa_n\Delta u_n=0.
\end{equation}
Again, $-\beta_n^{2n}\Delta \partial_t^{2n} u_n$ is a regularizing term.

\section{Heat Laws of Order $n$ with Memory}

\noindent
When $\kappa_n\neq 0$, we are in a position to make a memory relaxation of the law~\eqref{ordern},
along with a perturbation obtained by adding (minus) the gradient of the antiderivative $u_{n+1}$
of $u_n$. Hence, given a bounded convex summable function $g_{n+1}:\R^+\to\R^+$ of unitary total mass,
we consider the heat law \emph{of order $n$ with memory}
\begin{align}
\label{ordernmem}
\q+\sum_{i=1}^n  \alpha_n^i\partial_t^i\q
&= -\omega_{n+1}\kappa_n \nabla u_{n}
-(1-\omega_{n+1})\kappa_n \int_0^\infty g_{n+1}(s)\nabla u_{n}(t-s)ds\\
&\quad-\sum_{i=1}^{2n}\beta_n^i\nabla \partial_t^{i} u_{n}
-\kappa_{n+1}\nabla u_{n+1},\notag
\end{align}
where $\kappa_{n+1}\geq 0$ and $\omega_{n+1}\in[0,1)$.
This heat law, by applying the energy balance~\eqref{BASE}, yields the memory counterparts
of~\eqref{npiu}, exhibiting a nonlocal character in time due to the presence of the convolution integral.
Again, we shall distinguish two cases.

\medskip
\noindent
$\bullet$ If $\kappa_{n+1}=0$, we deduce 
the equation in the variable $u_n$
\begin{align}
\label{aaa}
&\sum_{i=1}^n  \alpha_n^i\partial_t^{i+n+1}u_{n}+\partial_t^{n+1}u_{n}
-\sum_{i=1}^{2n}\beta_n^i\Delta \partial_t^i u_n\\
&\quad-\omega_{n+1}\kappa_n \Delta u_{n}
-(1-\omega_{n+1})\kappa_n \int_0^\infty g_{n+1}(s)\Delta u_{n}(t-s)ds=0.\notag
\end{align}
As a particular instance of~\eqref{aaa} we consider the case
$n=1$ with $\beta_1^2=0$ (i.e., in absence of regularization), and all the other coefficients strictly positive, namely,
\begin{equation}
\label{tipoII}
\alpha_1^1\partial_{ttt}u_1+\partial_{tt}u_1-\beta_1^1\Delta\partial_t u_1
-\omega_{2}\kappa_1 \Delta u_{1}
-(1-\omega_{2})\kappa_1 \int_0^\infty g_{2}(s)\Delta u_{1}(t-s)ds=0.
\end{equation}
Introducing the integrated kernel
$$G(s)=(1-\omega_2)\kappa_1\int_s^\infty g_2(y)dy,$$
equation \eqref{tipoII} becomes
$$\alpha_1^1\partial_{ttt}u_1+\partial_{tt}u_1-\beta_1^1\Delta\partial_t u_1
-\kappa_1 \Delta u_{1}
+\int_0^\infty G(s)\Delta \partial_t u_{1}(t-s)ds=0,
$$
known in the literature as 
the MGT equation with memory of type II (see \cite{DLPII,LAW}).

\begin{remark}
As a matter of fact, we could also allow $\omega_2$ to take \emph{any} positive
value, except the forbidden value $\omega_2=1$, in compliance with Remark~\ref{remmoremmo}.
In particular, for $\omega_2>1$, equation~\eqref{tipoII}
becomes the MGT equation with memory of type I (see \cite{DLPI,LAW}).
\end{remark}

\medskip
\noindent
$\bullet$ If $\kappa_{n+1}\neq 0$, we write~\eqref{ordernmem}
in terms of $u_{n+1}$ only, and we deduce the equation
\begin{align}
\label{bbb}
&\sum_{i=1}^n  \alpha_n^i\partial_t^{i+n+2}u_{n+1}+\partial_t^{n+2}u_{n+1}
-\sum_{i=1}^{2n}\beta_n^i\Delta \partial_t^{i+1} u_{n+1}
-\omega_{n+1}\kappa_n \Delta \partial_t u_{n+1}\\
&\quad-(1-\omega_{n+1})\kappa_n \int_0^\infty g_{n+1}(s)\Delta \partial_t u_{n+1}(t-s)ds
-\kappa_{n+1}\Delta u_{n+1}=0.\notag
\end{align}
For $n=1$ and $\alpha_1^1=0$, hence $\beta_1^2=0$, equation~\eqref{bbb} reads
$$\partial_{ttt}u_2-\beta_1^1\Delta\partial_{tt} u_2
-\omega_{2}\kappa_1 \Delta \partial_t u_{2}
-(1-\omega_{2})\kappa_1 \int_0^\infty g_{2}(s)\Delta \partial_t u_{2}(t-s)ds
-\kappa_{2}\Delta u_{2}=0,
$$
which, defining the differentiated kernel
$$\mu(s)=-(1-\omega_{2})\kappa_1  g_{2}'(s),$$
can be written in the form
$$\partial_{ttt}u_2-\beta_1^1\Delta\partial_{tt} u_2
-\omega_{2}\kappa_1 \Delta \partial_t u_{2}
-\varrho\Delta u_{2}
+ \int_0^\infty \mu(s)\Delta u_{2}(t-s)ds,
$$
with
$\varrho=\kappa_2+(1-\omega_{2})\kappa_1g_2(0)$.
The latter is a regularized version of the 
MGT equation with memory of type I.

\section{The Inductive Procedure}

\noindent 
The last step is to justify the hierarchy of heat laws of order $n$ and of order $n$
with memory presented in the previous sections. This will be done by means of an inductive argument, starting from
the law of order $1$, whose derivation in turn starts from the Fourier law,
which in this scheme might be seen as the law of order $0$.
Accordingly, let us assume we are given for some $n\geq 1$ the law of order $n$ with memory~\eqref{ordernmem}.
We now choose the kernel $g_{n+1}$ of the form
$$
g(s)=\frac1{\varepsilon_{n+1}}\, e^{-\frac{s}{\varepsilon_{n+1}}},
$$
where the value $\varepsilon_{n+1}=0$ is allowed, and in that case $g_{n+1}$ becomes the Dirac mass
at zero.
At this point, exploiting~\eqref{dercon},
we take the time derivative of~\eqref{ordernmem},
written in terms of $u_{n+1}$ only, and we multiply it by 
$\varepsilon_{n+1}$.
This gives, changing the summation indexes,
\begin{align*}
\varepsilon_{n+1}\partial_t\q+\sum_{i=2}^{n+1}  \varepsilon_{n+1}\alpha_n^{i-1}\partial_t^{i}\q
&= -\varepsilon_{n+1}\omega_{n+1}\kappa_n \nabla \partial_{tt} u_{n+1}
-(1-\omega_{n+1})\kappa_n \nabla \partial_t u_{n+1}\\
&\quad +(1-\omega_{n+1})\kappa_n \int_0^\infty g_{n+1}(s)\nabla \partial_t u_{n+1}(t-s)ds\\
&\quad-\sum_{i=3}^{2n+2}\varepsilon_{n+1}\beta_n^{i-2}\nabla \partial_t^{i} u_{n+1}
-\varepsilon_{n+1}\kappa_{n+1}\nabla \partial_t u_{n+1}.
\end{align*}
Adding the latter equation with~\eqref{ordernmem} written
in terms of $u_{n+1}$, we finally get
(where the first sum disappears when $n=1$)
\begin{align*}
&\q+(\varepsilon_{n+1}+\alpha_n^1)\partial_t\q
+\sum_{i=2}^n  (\varepsilon_{n+1}\alpha_n^{i-1}+\alpha_n^i)\partial_t^i\q +\varepsilon_{n+1}\alpha_n^{n}\partial_t^{n+1}\q\\
&= -\kappa_{n+1}\nabla u_{n+1}-(\varepsilon_{n+1}\kappa_{n+1}+\kappa_n)\nabla\partial_t u_{n+1}
-(\varepsilon_{n+1}\omega_{n+1}\kappa_n+\beta_n^1)\nabla \partial_{tt} u_{n+1}\\
&\quad -\sum_{i=3}^{2n+1}(\varepsilon_{n+1}\beta_n^{i-2}+\beta_n^{i-1})\nabla \partial_t^{i} u_{n+1}-\varepsilon_{n+1}\beta_n^{2n}\nabla \partial_t^{2n+2} u_{n+1}.
\end{align*}
But, in the light of \eqref{alfa}-\eqref{beta}, 
this is exactly the law \eqref{ordern} written for $n+1$, that is,
the heat law of order $n+1$.
This completes the argument.

\section{Ten Equations}

\noindent
In this last section we produce a list of ten notable
evolution equations, all written in the generic variable $v$, that can be obtained from 
the family of constitutive laws for the heat flux introduced in this work.
In what follows, $a,b,c,d,e$ will stand for strictly positive constants, whereas
$g:\R^+\to\R^+$ will be a bounded convex summable function, not necessarily of unitary total mass.

\begin{itemize}

\medskip
\item[{\bf (i)}] {The heat equation:}
$$\partial_t v-a\Delta v=0.$$

\medskip
\item[{\bf (ii)}] {The Gurtin-Pipkin heat equation:}
$$\partial_t v-\int_0^\infty g(s)\Delta v(t-s)ds=0.$$

\medskip
\item[{\bf (iii)}] {The Coleman-Gurtin heat equation:}
$$\partial_t v-a\Delta v-\int_0^\infty g(s)\Delta v(t-s)ds=0.$$

\medskip
\item[{\bf (iv)}] {The weakly damped wave equation:}
$$\partial_{tt} v+a\partial_t v-b\Delta v=0.$$

\medskip
\item[{\bf (v)}] {The strongly damped wave equation:}
$$\partial_{tt} v-a \Delta \partial_t v-b\Delta v=0.$$

\medskip
\item[{\bf (vi)}] {The Moore-Gibson-Thompson (MGT) equation:}
$$\partial_{ttt} v+a\partial_{tt} v-b\Delta \partial_{t} v -c\Delta v=0.$$

\medskip
\item[{\bf (vii)}] {The regularized MGT equation:}
$$\partial_{ttt} v+a\partial_{tt} v-b\Delta \partial_{tt} v -c\Delta \partial_{t}v
-d\Delta v=0.$$

\medskip
\item[{\bf (viii)}] {The MGT equation with memory of type I:}
$$\partial_{ttt} v+a\partial_{tt} v-b\Delta \partial_{t} v -c\Delta v
+\int_0^\infty g(s)\Delta v(t-s)ds=0.$$

\medskip
\item[{\bf (ix)}] {The MGT equation with memory of type II:}
$$\partial_{ttt} v+a\partial_{tt} v-b\Delta \partial_{t} v -c\Delta v
+\int_0^\infty g(s)\Delta \partial_t v(t-s)ds=0.$$

\medskip
\item[{\bf (x)}] {A fourth order equation of MGT type:}
$$\partial_{tttt} v+a\partial_{ttt} v+b\partial_{tt} v
-c \Delta \partial_{tt} v -d \Delta \partial_{t} v
-e\Delta v=0.$$

\end{itemize}



\end{document}